\documentclass{fundam}

\usepackage{cite}
\usepackage{mathrsfs}
\usepackage{amsfonts}
\usepackage{amsmath}
\usepackage{amsfonts,amssymb}
\usepackage{dsfont}
\usepackage{mathrsfs}
\usepackage{pifont}
\usepackage{amssymb}
\usepackage{graphicx}
\allowdisplaybreaks

\numberwithin{equation}{section}

\usepackage{hyperref}

\begin{document}

%%%%%%%  parameters to be filled in by copy-editor  %%%%%%%%%%

\setcounter{page}{67}
\publyear{24}
\papernumber{2172}
\volume{191}
\issue{1}

\finalVersionForARXIV
%\finalVersionForIOS

%%%%%%%%%%%%%%%%%%%%%%%%%%%%%%%%%%%%%%

\title{Two Sufficient Conditions for Graphs to Admit Path Factors}

\author{Sizhong Zhou, Jiancheng Wu\thanks{Address for correspondence:  School of Science, Jiangsu
            University of Science and Technology,  Zhenjiang, Jiangsu 212100, China. \newline \newline
             \vspace*{-6mm}{\scriptsize{Received August 2023; \ accepted January 2024.}}}
\\
 School of Science\\
 Jiangsu University of Science and Technology\\
 Zhenjiang, Jiangsu 212100, China\\
 zsz\_cumt@163.com, wujiancheng@just.edu.cn
  }

\maketitle
\runninghead{S. Zhou and J. Wu}{Two Sufficient Conditions for Graphs to Admit Path Factors}

\vspace*{-4mm}
\begin{abstract}
\noindent Let $\mathcal{A}$ be a set of connected graphs. Then a spanning subgraph $A$ of $G$ is called an $\mathcal{A}$-factor if each component
of $A$ is isomorphic to some member of $\mathcal{A}$. Especially, when every graph in $\mathcal{A}$ is a path, $A$ is a path factor. For a positive
integer $d\geq2$, we write $\mathcal{P}_{\geq d}=\{P_i|i\geq d\}$. Then a $\mathcal{P}_{\geq d}$-factor means a path factor in which every component
admits at least $d$ vertices. A graph $G$ is called a $(\mathcal{P}_{\geq d},m)$-factor deleted graph if $G-E'$ admits a $\mathcal{P}_{\geq d}$-factor
for any $E'\subseteq E(G)$ with $|E'|=m$. A graph $G$ is called a $(\mathcal{P}_{\geq d},k)$-factor critical graph if $G-Q$ has a
$\mathcal{P}_{\geq d}$-factor for any $Q\subseteq V(G)$ with $|Q|=k$. In this paper, we present two degree conditions for graphs to be
$(\mathcal{P}_{\geq3},m)$-factor deleted graphs and $(\mathcal{P}_{\geq3},k)$-factor critical graphs. Furthermore, we show that the two results
are best possible in some sense.
\end{abstract}


\section{Introduction}

In our daily life many physical structures can conveniently be simulated by networks. The core issue of network security is the ruggedness and
vulnerability of the network, which is also one of the key topics that researchers must consider during the network designing phase. To study
the properties of the network, we use a graph to simulate the network. Vertices of the graph stand for nodes of the network and edges of the
graph act as links between the nodes of the network. Henceforth, we replace ``network" by the term ``graph". In data transmission networks,
the data transmission between two nodes of a network stands for a path between two corresponding vertices of a corresponding graph. Consequently,
the availability of data transmission in the network is equivalent to the existence of path factor in the corresponding graph which is generated
by the network. Clearly, research on the existence of path factors under specific network structures can help scientists design and construct
networks with high data transmission rates. Furthermore, the existence of a path factor deleted graph and a path factor critical graph
also plays a key role in data transmission of a network. If some links (resp. nodes) are damaged in the process of data transmission at the moment,
the possibility of data transmission between nodes is characterized by whether the corresponding graph of the network is a path factor deleted
(resp. critical) graph. In this article, we investigate the existence of path factor deleted graphs and path factor critical graphs which play a
key role in studying data transmissions of data transmission networks. We find that there are strong essential connection between some graphic
parameters (for instance, degree and connectivity, and so on) and the existence of path factor deleted graphs (or path factor critical graphs),
and hence investigations on degree and connectivity, which play an irreplaceable role in the vulnerability of the network and the feasibility of
data transmission, can yield theoretical guidance to meet data transmission and network security requirements.

\medskip
We discuss only finite undirected graphs without loops or multiple edges, otherwise, unless expli\-citly stated. We use $G\!=\!(V(G),E(G))$ to denote a graph, where $V(G)$ is the vertex set of $G$ and $E(G)$ is the edge set of $G$. The degree of a vertex $x$ in $G$, denoted by $d_G(x)$, is the number
of edges incident with $x$ in $G$. The neighborhood of a vertex $x$ in $G$, denoted by $N_G(x)$, is the set of vertices adjacent to $x$ in $G$.
For any $X\subseteq V(G)$, we write $N_G(X)=\bigcup\limits_{x\in X}{N_{G}(x)}$ and denote by $G[X]$ the subgraph of $G$ induced by $X$. Let
$G-X=G[V(G)\setminus X]$. For an edge subset $E'$ of $G$, we use $G-E'$ to denote the subgraph obtained from $G$ by removing $E'$. If $d_G(x)=0$
for some vertex $x$ in $G$, then $x$ is called an isolated vertex in $G$. We denote by $I(G)$ the set of isolated vertices of $G$, and set
$i(G)=|I(G)|$. The number of connected components of $G$ is denoted by $\omega(G)$. We denote by~$\lambda(G)$ and $\kappa(G)$ the edge connectivity and the vertex connectivity of $G$, respectively. A vertex subset $X$ of $G$ is said to be independent if $G[X]$ has no edges. The binding number of $G$ is defined by Woodall \cite{W}~as
$$
bind(G)=\min\left\{\frac{|N_G(X)|}{|X|}:\emptyset\neq X\subseteq V(G),N_G(X)\neq V(G)\right\}.
$$
We denote by $P_n$ and $K_n$ the path and the complete graph of order $n$, respectively. Let $G_1$ and $G_2$ be two graphs. Then we denote by
$G_1\cup G_2$ and $G_1\vee G_2$ the union and the join of $G_1$ and $G_2$, respectively.

\medskip
Let $\mathcal{A}$ be a set of connected graphs. Then a spanning subgraph $A$ of $G$ is called an $\mathcal{A}$-factor if each component of $A$
is isomorphic to some member of $\mathcal{A}$. Especially, when every graph in $\mathcal{A}$ is a path, $A$ is a path factor. For a positive
integer $d\geq2$, we write $\mathcal{P}_{\geq d}=\{P_i|i\geq d\}$. Then a $\mathcal{P}_{\geq d}$-factor means a path factor in which every
component admits at least $d$ vertices. In order to characterize a graph with a $\mathcal{P}_{\geq3}$-factor, Kaneko \cite{K} introduced the
concept of a sun. A graph $H$ is called a factor-critical graph if every induced subgraph with $|V(H)|-1$ vertices of $H$ admits a perfect
matching. Assume that $H$ is a factor-critical graph with vertex set $V(H)=\{v_1,v_2,\cdots,v_n\}$. By adding new vertices $u_1,u_2,\cdots,u_n$
and new edges $v_1u_1,v_2u_2,\cdots,v_nu_n$ to $H$, we acquire a new graph $R$, which is called a sun. In view of Kaneko, $K_1$ and $K_2$ are
also suns. A sun with at least six vertices is said to be a big sun. We denote by $sun(G)$ the number of sun components of $G$.

%%\medskip
Las Vergnas \cite{LV} characterized a graph admitting a $\mathcal{P}_{\geq2}$-factor. Bazgan, Benhamdine, Li and Wo\'niak \cite{BBLW} claimed
that a 1-tough graph $G$ of order at least 3 admits a $\mathcal{P}_{\geq3}$-factor. Kaneko \cite{K} derived a criterion for a graph with a
$\mathcal{P}_{\geq3}$-factor. Kano, Katona and Kir\'aly \cite{KKK} gave a simpler proof of Kaneko's result. Kano, Lu and Yu \cite{KLY}
demonstrated that a graph $G$ has a $\mathcal{P}_{\geq3}$-factor if $i(G-X)\leq\frac{2}{3}|X|$ for all $X\subseteq V(G)$. Liu \cite{L} showed
a result on the existence of $\mathcal{P}_{\geq3}$-factors in graphs with given properties. Wang and Zhang \cite{WZi}, Wu \cite{Wp}, Zhou et al
\cite{Zs,ZSL1,ZZS,ZSB,ZSY,ZWB} posed some sufficient conditions for graphs having $\mathcal{P}_{\geq3}$-factors. Gao, Wang and Chen \cite{GWC}
obtained some tight bounds for the existence of path factors in graphs. Kano, Lee and Suzuki \cite{KLS} verified that a connected cubic graph
with at least eight vertices admits a $\mathcal{P}_{\geq8}$-factor. Wang and Zhang \cite{WZs}, Katerinis and Woodall \cite{KW}, Zhou, Bian and
Sun \cite{ZBS} established some relationships between binding numbers and graph factors. Wu \cite{Wa}, Wang and Zhang \cite{WZd,WZo}, Zhou, Pan
and Xu \cite{ZPX}, Zhou \cite{Za,Zd}, Zhou, Wu and Liu \cite{ZWL} presented some degree conditions for the existence of graph factors in graphs.
Some other results on graph factors can be seen at \cite{WZ,ZSL2,ZL2}.

The following theorem on $\mathcal{P}_{\geq3}$-factors are known, which plays an important role in the proofs of our main results.

\medskip
\noindent{\textbf{Theorem 1}} (\cite{K,KKK}). A graph $G$ has a $\mathcal{P}_{\geq3}$-factor if and only if
$$
sun(G-X)\leq2|X|
$$
for any vertex subset $X$ of $G$.

\medskip
A graph $G$ is called a $(\mathcal{P}_{\geq d},m)$-factor deleted graph if $G-E'$ admits a $\mathcal{P}_{\geq d}$-factor for any $E'\subseteq E(G)$
with $|E'|=m$. A graph $G$ is called a $(\mathcal{P}_{\geq d},k)$-factor critical graph if $G-Q$ has a $\mathcal{P}_{\geq d}$-factor for any
$Q\subseteq V(G)$ with $|Q|=k$. Zhou \cite{Zr} showed two binding number conditions for graphs to be $(\mathcal{P}_{\geq3},m)$-factor deleted graphs
and $(\mathcal{P}_{\geq3},k)$-factor critical graphs.

\medskip
\noindent{\textbf{Theorem 2}} (\cite{Zr}). Let $m$ be a nonnegative integer, and let $G$ be a graph. If $\kappa(G)\geq2m+1$ and
$bind(G)>\frac{3}{2}-\frac{1}{4m+4}$, then $G$ is a $(\mathcal{P}_{\geq3},m)$-factor deleted graph.

\medskip
\noindent{\textbf{Theorem 3}} (\cite{Zr}). Let $k$ be a nonnegative integer, and let $G$ be a graph with $\kappa(G)\geq k+2$. If $bind(G)\geq\frac{5+k}{4}$,
then $G$ is a $(\mathcal{P}_{\geq3},k)$-factor critical graph.

\medskip
It is natural and interesting to put forward some new graphic parameter conditions to ensure that a graph is a $(\mathcal{P}_{\geq3},m)$-factor
deleted graph or a $(\mathcal{P}_{\geq3},k)$-factor critical graph. In this paper, we show degree conditions for graphs to be $(\mathcal{P}_{\geq3},m)$-factor
deleted graphs or $(\mathcal{P}_{\geq3},k)$-factor critical graphs, which are shown in the following.

\medskip
\noindent{\textbf{Theorem 4.}} Let $m$ and $r$ be two integers with $r\geq1$ and $0\leq m\leq2r+1$, and let $G$ be a graph of order $n$ with
$n\geq4r+6m+4$. If $\kappa(G)\geq r+m$ and $G$ satisfies
$$
\max\{d_G(x_1),d_G(x_2),\cdots,d_G(x_{2r+1})\}\geq\frac{n}{3}
$$
for any independent set $\{x_1,x_2,\cdots,x_{2r+1}\}$ of $G$, then $G$ is a $(\mathcal{P}_{\geq3},m)$-factor deleted graph.

\eject
%%%\medskip
\noindent{\textbf{Theorem 5.}} Let $k\geq0$ and $r\geq1$ be two integers, and let $G$ be a graph of order $n$ with $n\geq4r+k+4$. If
$\kappa(G)\geq r+k$ and $G$ satisfies
$$
\max\{d_G(x_1),d_G(x_2),\cdots,d_G(x_{2r+1})\}\geq\frac{n+2k}{3}
$$
for any independent set $\{x_1,x_2,\cdots,x_{2r+1}\}$ of $G$, then $G$ is a $(\mathcal{P}_{\geq3},k)$-factor critical graph.

\medskip

From a computation standpoint, the condition given in each of the two theorems can be checked in polynomial time for fixed $r$, $m$ and $k$.

\section{The proof of Theorem 4}

\noindent{\it Proof of Theorem 4.} Let $G'=G-E'$ for $E'\subseteq E(G)$ with $|E'|=m$. Then $V(G')=V(G)$ and $E(G')=E(G)\setminus E'$. To prove
Theorem 4, it suffices to verify that $G'$ admits a $\mathcal{P}_{\geq3}$-factor. Suppose, to the contrary, that $G'$ has no $\mathcal{P}_{\geq3}$-factor.
Then it follows from Theorem 1 that
$$
sun(G'-X)\geq2|X|+1\eqno(1)
$$
for some $X\subseteq V(G')$.

\medskip
Next, we shall consider two cases according to the value of $i(G-X)$ and derive a contradiction in each case.

\medskip
\noindent{\bf Case 1.} $i(G-X)\geq2r+1$.

\smallskip
In this case, $G-X$ admits at least $2r+1$ isolated vertices $v_1,v_2,\cdots,v_{2r+1}$. Thus, we acquire $d_{G-X}(v_i)=0$ for $1\leq i\leq2r+1$,
and so
$$
d_G(v_i)\leq d_{G-X}(v_i)+|X|=|X|\eqno(2)
$$
for $1\leq i\leq2r+1$.

\medskip
Obviously, $\{v_1,v_2,\cdots,v_{2r+1}\}$ is an independent set of $G$. Then using (2) and the degree condition of Theorem 4, we infer
$$
|X|\geq\max\{d_G(v_1),d_G(v_2),\cdots,d_G(v_{2r+1})\}\geq\frac{n}{3}.\eqno(3)
$$

It follows from (1) and (3) that
$$
n\geq|X|+sun(G'-X)\geq|X|+2|X|+1=3|X|+1\geq3\cdot\frac{n}{3}+1=n+1,
$$
which is a contradiction.

\medskip
\noindent{\bf Case 2.} $i(G-X)\leq2r$.

\medskip
\noindent{\bf Subcase 2.1.} $X$ is not a vertex cut set of $G$.

\medskip
Clearly, $\omega(G-X)=\omega(G)=1$. If $|X|\geq\frac{m+1}{2}$, then we have
$$
sun(G'-X)=sun(G-X-E')\leq\omega(G-X-E')\leq\omega(G-X)+m=m+1\leq2|X|,
$$
which contradicts (1).
\eject

%%%\medskip
If $1\leq|X|<\frac{m+1}{2}$, then it follows from $m\leq2r+1$ that
\begin{eqnarray*}
            \lambda(G-X) & \geq& \kappa(G-X)\geq\kappa(G)-|X| >  (r+m)-\frac{m+1}{2}   \nonumber \\
            & \geq & \left(\frac{m-1}{2}+m\right)-\frac{m+1}{2}=m-1.
\end{eqnarray*}
By the integrity of $\lambda(G-X)$, we get
$$
\lambda(G-X)\geq m.\eqno(4)
$$

For $e\in E'$, we admit by (4)
$$
\lambda(G'-X+\{e\})=\lambda(G-X-E'\setminus\{e\})\geq\lambda(G-X)-(m-1)\geq m-(m-1)=1,
$$
which implies $\omega(G'-X+\{e\})=1$. Hence, we derive
$$
\omega(G'-X)\leq\omega(G'-X+\{e\})+1=2.\eqno(5)
$$

Using (1), (5) and $1\leq|X|<\frac{m+1}{2}$, we deduce
$$
3\leq2|X|+1\leq sun(G'-X)\leq\omega(G'-X)\leq2,
$$
which is a contradiction.

\medskip
If $|X|=0$, then from (1), we have
$$
sun(G')=sun(G'-X)\geq2|X|+1=1.\eqno(6)
$$
Note that $\kappa(G)\geq r+m$ and $|E'|=m$. Then $G'=G-E'$ is $r$-connected, and so $1=\omega(G')\geq sun(G')$. Combining this with (6), we get
$$
sun(G')=1.\eqno(7)
$$

In view of (7) and $n\geq4r+6m+4$, we know that $G'$ is a big sun. Let $R$ denote the factor-critical graph of $G'$ with $|V(R)|=\frac{n}{2}$.
Then $G'$ has at least $2r+3m+2$ vertices with degree 1. Note that $G'=G-E'$ with $|E'|=m$. Hence, $G$ admits at least $2r+m+2$ vertices with
degree 1. Then we may select an independent set $\{v_1,v_2,\cdots,v_{2r+1}\}\subseteq V(G)\setminus V(R)$ such that $d_G(v_i)=1$ for
$1\leq i\leq2r+1$. According to the degree condition of Theorem 4, we infer
$$
1=\max\{d_G(v_1),d_G(v_2),\cdots,d_G(v_{2r+1})\}\geq\frac{n}{3},
$$
that is,
$$
n\leq3,
$$
which contradicts $n\geq4r+6m+4$.

\eject
\noindent \textbf{Subcase 2.2.} $X$ is a vertex cut set of $G$.

\medskip
In this subcase, $\omega(G-X)\geq2$ and $|X|\geq r+m$. In terms of (1), we obtain
\begin{eqnarray*}
 sun(G-X) \hspace*{-6mm} && \geq  sun(G-X-E')-2m   \nonumber  \\
& & =  sun(G'-X)-2m\geq2|X|+1-2m\geq2(r+m)+1-2m   \nonumber  \\
& & = 2r+1,
\end{eqnarray*}
which implies that there exist $t$ sun components in $G-X$, denoted by $H_1,H_2,\cdots,H_t$, where $t\geq2r+1$. We choose $v_i\in V(H_i)$ with
$d_{H_i}(v_i)\leq1$ for $1\leq i\leq2r+1$. Clearly, $\{v_1,v_2,\cdots,v_{2r+1}\}$ is an independent set of $G$. Then it follows from the degree
condition of Theorem 4 that
$$
\max\{d_G(v_1),d_G(v_2),\cdots,d_G(v_{2r+1})\}\geq\frac{n}{3}.\eqno(8)
$$
Without loss of generality, we may let $d_G(v_1)\geq\frac{n}{3}$ by (8). Thus, we infer
$$
d_{G[X]}(v_1)=d_G(v_1)-d_{H_1}(v_1)\geq\frac{n}{3}-1,
$$
and so
$$
|X|\geq d_{G[X]}(v_1)\geq\frac{n}{3}-1.\eqno(9)
$$

In light of (1), (9), $r\geq1$, $i(G-X)\leq2r$ and $n\geq4r+6m+4$, we deduce
\begin{eqnarray*}
n&\geq&|X|+2\cdot sun(G-X)-i(G-X)\\[2pt]
&\geq&|X|+2(sun(G-X-E')-2m)-i(G-X)\\[2pt]
&=&|X|+2(sun(G'-X)-2m)-i(G-X)\\[2pt]
&\geq&|X|+2(2|X|+1-2m)-2r\\[2pt]
&=&5|X|-4m-2r+2\\[2pt]
&\geq&5\left(\frac{n}{3}-1\right)-4m-2r+2\\[2pt]
&=&n+\frac{2n}{3}-4m-2r-3\\[2pt]
&\geq&n+\frac{2(4r+6m+4)}{3}-4m-2r-3\\[2pt]
&=&n+\frac{2r}{3}-\frac{1}{3}\\[2pt]
&\geq&n+\frac{1}{3}\\[2pt]
&>&n,
\end{eqnarray*}

\noindent which is a contradiction. We complete the proof of Theorem 4. \hfill $\Box$

\medskip
\noindent{\bf Remark 1.} We now show that
$$
\max\{d_G(x_1),d_G(x_2),\cdots,d_G(x_{2r+1})\}\geq\frac{n}{3}
$$
\eject
\noindent in Theorem 4 cannot be replaced by
$$
\max\{d_G(x_1),d_G(x_2),\cdots,d_G(x_{2r+1})\}\geq\frac{n-1}{3}.
$$

Let $m\geq0$ and $r\geq1$ be two integers, and $t$ be a sufficiently large integer. We construct a graph $G=K_{rt+m}\vee((2rt+1)K_1\cup (mK_2))$.
Then $G$ is $(rt+m)$-connected, $n=3rt+3m+1$ and
$$
\max\{d_G(x_1),d_G(x_2),\cdots,d_G(x_{2r+1})\}\geq rt+m=\frac{n-1}{3}
$$
for any independent set $\{x_1,x_2,\cdots,x_{2r+1}\}$ of $G$. Let $E'=E(mK_2)$ and $G'=G-E'$. Choose $X=V(K_{rt+m})$. Then we obtain
$$
sun(G'-X)=2rt+2m+1=2|X|+1>2|X|.
$$
By virtue of Theorem 1, $G'$ has no $\mathcal{P}_{\geq3}$-factor, and so $G$ is not a $(\mathcal{P}_{\geq3},m)$-factor deleted graph.

\medskip
But we do not know whether the condition $\kappa(G)\geq r+m$ in Theorem 4 is best possible or not. Naturally, it is interesting to further study
the above problem.

\section{The proof of Theorem 5}

\noindent{\it Proof of Theorem 5.} Let $G'=G-Q$ for $Q\subseteq V(G)$ with $|Q|=k$. To verify Theorem 5, it suffices to claim that $G'$ has a
$\mathcal{P}_{\geq3}$-factor. Suppose, to the contrary, that $G'$ has no $\mathcal{P}_{\geq3}$-factor. Then by Theorem~1, we derive
$$
sun(G'-X)\geq2|X|+1\eqno(1)
$$
for some $X\subseteq V(G')$.

\medskip
\noindent \textbf{Claim 1.} $|X|\geq r$.

\begin{proof}
 Assume $|X|\leq r-1$. Since $\kappa(G)\geq r+k$ and $|Q|=k$, $G'-X=G-Q-X$ is connected, which implies $\omega(G'-X)=1$.
By (1), we get
$$
1\leq2|X|+1\leq sun(G'-X)\leq\omega(G'-X)=1,
$$
which implies $|X|=0$ and $sun(G'-X)=sun(G')=1$. Then by $|V(G')|=n-k\geq4r+4$, we see that $G'$ is a big sun. Let $R$ be the factor-critical
graph of $G'$ with $|V(R)|=\frac{1}{2}|V(G')|=\frac{1}{2}(n-k)$. Then by $n\geq4r+k+4$, $G'$ admits an independent set
$\{v_1,v_2,\cdots,v_{2r+1}\}\subseteq V(G')\setminus V(R)$ with $d_{G'}(v_i)=1$ for $1\leq i\leq 2r+1$.

\smallskip According to the degree condition of Theorem 5, we obtain
\begin{eqnarray*}
k+1&=&|Q|+1=|Q|+\max\{d_{G'}(v_1),d_{G'}(v_2),\cdots,d_{G'}(v_{2r+1})\}\\
&\geq&\max\{d_G(v_1),d_G(v_2),\cdots,d_G(v_{2r+1})\}\\
&\geq&\frac{n+2k}{3},
\end{eqnarray*}
\eject

\noindent that is,
$$
n\leq k+3,
$$
which contradicts $n\geq4r+k+4$. Claim 1 is proved.
\end{proof}

\noindent \textbf{Claim 2.} $i(G'-X)\leq2r$.

\begin{proof}
Let $i(G'-X)\geq2r+1$. Then there exist at least $2r+1$ isolated vertices $v_1,v_2,\cdots,v_{2r+1}$
in $G'-X$, namely, $d_{G'-X}(v_i)=0$ for $1\leq i\leq2r+1$. Hence, we have
$$
d_G(v_i)\leq d_{G-Q}(v_i)+|Q|=d_{G'}(v_i)+k\leq d_{G'-X}(v_i)+|X|+k=|X|+k\eqno(2)
$$
for $1\leq i\leq2r+1$.

\medskip
In terms of (2) and the degree condition of Theorem 5, we obtain
$$
|X|+k\geq\max\{d_G(v_1),d_G(v_2),\cdots,d_G(v_{2r+1})\}\geq\frac{n+2k}{3},
$$
which implies
$$
|X|\geq\frac{n-k}{3}.\eqno(3)
$$

It follows from (1) and (3) that
\begin{eqnarray*}
   n &\geq&  |Q|+|X|+sun(G'-X)\geq k+|X|+2|X|+1\nonumber\\
  & = & 3|X|+k+1\geq3\cdot\frac{n-k}{3}+k+1=n+1,
\end{eqnarray*}
which is a contradiction. We complete the proof of Claim 2.
\end{proof}

Using (1) and Claim 1, we derive
$$
sun(G'-X)\geq2|X|+1\geq2r+1,
$$
which implies that $G'-X$ admits $t$ sun components $H_1,H_2,\cdots,H_t$, where $t\geq2r+1$. Choose $v_i\in V(H_i)$ such that $d_{H_i}(v_i)\leq1$
for $1\leq i\leq2r+1$. Obviously, $\{v_1,v_2,\cdots,v_{2r+1}\}$ is an independent set of $G$. By the degree condition of Theorem 5, we have
$$
\max\{d_G(v_1),d_G(v_2),\cdots,d_G(v_{2r+1})\}\geq\frac{n+2k}{3}.\eqno(4)
$$
Without loss of generality, let $d_G(v_1)\geq\frac{n+2k}{3}$ by (4). Thus, we derive
$$
k+|X|=|Q|+|X|\geq d_{G[Q\cup X]}(v_1)=d_G(v_1)-d_{H_1}(v_1)\geq\frac{n+2k}{3}-1,
$$
that is,
$$
|X|\geq\frac{n-k}{3}-1.\eqno(5)
$$

According to (1), (5), Claim 2, $r\geq1$ and $n\geq4r+k+4$, we obtain
\begin{eqnarray*}
n&\geq&|Q|+|X|+2\cdot sun(G'-X)-i(G'-X)\\
&\geq&k+|X|+2(2|X|+1)-2r\\
&=&5|X|+k-2r+2\\
&\geq&5\left(\frac{n-k}{3}-1\right)+k-2r+2\\
&=&n+\frac{2n}{3}-\frac{2k}{3}-2r-3\\
&\geq&n+\frac{2(4r+k+4)}{3}-\frac{2k}{3}-2r-3\\
&=&n+\frac{2r}{3}-\frac{1}{3}\\
&\geq&n+\frac{1}{3}\\
&>&n,
\end{eqnarray*}
which is a contradiction. The proof of Theorem 5 is complete. \hfill $\Box$

\medskip

\noindent{\bf Remark 2.} We now show that
$$
\max\{d_G(x_1),d_G(x_2),\cdots,d_G(x_{2r+1})\}\geq\frac{n+2k}{3}
$$
in Theorem 5 cannot be replaced by
$$
\max\{d_G(x_1),d_G(x_2),\cdots,d_G(x_{2r+1})\}\geq\frac{n+2k-1}{3}.
$$

Let $k\geq0$ and $r\geq1$ be two integers, and $t$ be a sufficiently large integer. We construct a graph $G=K_{rt+2k+1}\vee((2rt+2k+3)K_1)$.
Then $G$ is $(rt+2k+1)$-connected, $n=3rt+4k+4=3(rt+2k+1)-2k+1$ and
$$
\max\{d_G(x_1),d_G(x_2),\cdots,d_G(x_{2r+1})\}=rt+2k+1=\frac{n+2k-1}{3}
$$
for any independent subset $\{x_1,x_2,\cdots,x_{2r+1}\}$ of $G$. Let $Q\subseteq V(K_{rt+2k+1})$ with $|Q|=k$, and $G'=G-Q$. Select
$X=V(K_{rt+k+1})\subseteq V(K_{rt+2k+1})\setminus Q$. Then we possess
$$
sun(G'-X)=2rt+2k+3=2(rt+k+1)+1=2|X|+1>2|X|.
$$
Using Theorem 1, $G'$ has no $\mathcal{P}_{\geq3}$-factor, and so $G$ is not a $(\mathcal{P}_{\geq3},k)$-factor critical graph.

\medskip
But we do not know whether the condition $\kappa(G)\geq r+k$ in Theorem 5 is best possible or not. Naturally, it is interesting to further study
the above problem.

\section*{Data availability statement}

My manuscript has no associated data.

\section*{Declaration of competing interest}

The authors declare that they have no conflicts of interest to this work.

\subsection*{Acknowledgments}

The authors would like to express their gratitude to the anonymous reviewers for their helpful
comments and valuable suggestions in improving this paper. This work was supported by the
Natural Science Foundation of Shandong Province, China (ZR2023MA078).

\end{document}